\documentclass[10pt, reqno]{amsart}
\usepackage{fancyhdr}
\usepackage{a4wide, amsthm, amscd, amsfonts, amssymb, graphicx,tikz, color, environ}
\usepackage{amsaddr}
\usepackage{amsthm}

\usepackage{amsmath,amssymb,amsfonts,graphicx}
\usepackage{amsmath,amssymb,amsfonts,amsthm, graphicx,color,fancyhdr}
\usepackage[latin1]{inputenc}
\usepackage{enumerate}

\usepackage{psfrag}
\usepackage[colorlinks=true, pdfstartview=FitV, allcolors=black]{hyperref}
\usepackage{color}
\usepackage{psfrag,caption}
\usepackage{amsmath}
\usepackage{amsthm}
\usepackage{amscd,amssymb}
\usepackage[all]{xy}
\usepackage{color}
\usepackage{comment}
\usepackage{appendix}
\usepackage{verbatim}
\usepackage{mathrsfs}
\usepackage{pdfsync}
\usepackage{graphicx,epsfig}
\newtheorem{thm}{Theorem}[section]

\theoremstyle{stylename}

\let\oldproofname=\proofname
\renewcommand{\proofname}{\rm\bf{\oldproofname}}

 \newtheorem{lemma}[thm]{Lemma}

\theoremstyle{definition}
 
 \newtheorem{rmk}[thm]{Remark}
 \newtheorem{defn}{Definition}[section]

\subjclass[2010]{35P15, 53C20, 53C42, 58J50}

\begin{document}
\title{Bounds for the first non-zero Steklov eigenvalue}
\author{Sheela Verma}
\address{Department of Mathematics and Statistics\\ Indian Institute of Technology Kanpur\\ Kanpur, India}
\email{sheela@iitk.ac.in}

\begin{abstract}
Let $\Omega$ be a star-shaped bounded domain in $(\mathbb{S}^{n}, ds^{2})$ with smooth boundary. In this article, we give a sharp lower bound for the first non-zero eigenvalue of the Steklov eigenvalue problem in $\Omega.$ This result is the generalization of a result given by Kuttler and Sigillito for a star-shaped bounded domain in $\mathbb{R}^2.$ Further we also obtain a two sided bound for the first non-zero eigenvalue of the Steklov problem on the ball in $\mathbb{R}^n$ with rotationally invariant metric and with bounded radial curvature.
\end{abstract}

\keywords {Laplacian; Steklov eigenvalue problem; Radial curvature; Rotationally invariant metric. }
\maketitle

\section{Introduction}
Let $\Omega $ be a bounded domain in $(\mathbb{S}^{n}, ds^{2})$ with smooth boundary $\partial\Omega$. Consider the following problem
\begin{align} \label{Steklov problem}
\begin{array}{rcll}
\Delta f &=& 0 & \mbox{ in } \Omega ,\\
\frac{\partial f}{\partial \nu} &=& \mu f  &\mbox{ on } \partial \Omega,
\end{array}
\end{align}
where $\nu$ is the unit outward normal on the boundary $\partial \Omega$ and $\mu$ is a real number.

This problem is known as Steklov eigenvalue problem and was introduced by Steklov \cite{S} for bounded domains in the plane in $1902.$ This problem is important as the set of eigenvalues of the Steklov problem is same as the set of eigenvalues of the well known Dirichlet-Neumann map. This map associates to each function $f$ defined on $\partial \Omega$, the normal derivative of its harmonic extension on $\Omega$. The eigenvalues of the Steklov problem are discrete and forms an increasing sequence $0= \mu_0 < \mu_1 < \mu_2 < \cdots$.

In this article, we are interested in finding bounds for the first non-zero eigenvalue $\mu_1$ of the Steklov problem. The variational characterization of $\mu_1$ is given by
\begin{align} \label{variational characterization}
\mu_{1}= \inf \left\lbrace \frac{\int_\Omega{\|\nabla u\|^2}\, dv}{\int_{\partial\Omega}{u^2}\, ds} : \int_{\partial\Omega}{u}\, ds =0, \, u(\neq 0) \in C^1(\overline{\Omega})\right\rbrace.
\end{align}

There are several results which estimate the first non-zero eigenvalue $\mu_1$ of the Steklov eigenvalue problem \cite{BP, E1, E2, E3, BG}. 
The first upper bound for $\mu_{1}$ was given by Weinstock \cite{W} in $1954$. He proved that among all simply connected planar domains with analytic boundary of fixed perimeter, the circle maximizes $\mu_{1}$.

Let $\Omega\subset \mathbb{R}^{n}$ be a star shaped domain with smooth boundary $\partial\Omega$. Let $p$ be a centre of $\Omega$. Let $R_m = \min\left\lbrace  d(p, x) | x \in \partial\Omega \right\rbrace $, $R_M = \max\left\lbrace d(p, x) | x \in \partial\Omega\right\rbrace $ and $h_m = \min\left\lbrace \langle x, \nu \rangle | x \in \partial\Omega \right\rbrace $, where $\nu$ is the outward unit normal to $\partial\Omega$. With these notations Bramble and Payne \cite{BP} proved that

$$ \mu_{1} \geq \frac{{R_m}^{n-1}}{{R_M}^{n+1}} \, h_m. $$
This bound is sharp when $\Omega$ is a ball. Payne \cite{P} calculated the two sided bound for $\mu_1$ on convex domain $\Omega \subset \mathbb{R}^{2}$ in terms of minimum and maximum curvature on boundary of $\Omega$. Escobar \cite{E1} generalized this result for two dimensional compact Riemannian manifold $(M,g)$ with non-negative Gaussian curvature such that geodesic curvature on $\partial M$ is bounded from below.

With a different approach, Kuttler and Sigillito \cite{KS} proved the following theorem for a star-shaped bounded domain in $\mathbb{R}^2$.

\begin{thm}[\cite{KS}]
Let $\Omega$ be a star-shaped bounded domain in $\mathbb{R}^2$ with smooth boundary. Then 
\begin{align*}
\mu_1(\Omega) \geq \dfrac{\left[ 1- \dfrac{2}{1+\sqrt{1+ 4\, \min \left( \frac{R(\theta)}{R'(\theta)}\right)^2}}\right] }{\max \sqrt{R^2(\theta) + {R'}^2(\theta)}},  
\end{align*}
where $R(\theta):= \max \left\lbrace |x| : x \in \Omega, x = |x| e^{i\theta} \right\rbrace $  and equality holds for a disc. 
\end{thm}

Garcia and Montano \cite{GM} generalized this result for the star-shaped domain in $\mathbb{R}^n$. 
\begin{thm}[\cite{GM}]
Let $\Omega \subset \mathbb{R}^n$ be a star-shaped domain with smooth boundary $\partial\Omega$ and unit outward normal $\nu.$ If $0 \leq \theta \leq \alpha < \frac{\pi}{2}$, where $\cos (\theta) = \langle \nu, \partial_{r}\rangle$, then the first non-zero eigenvalue of the Steklov problem $\mu_{1}(\Omega)$ satisfies 
\begin{align*}
\mu_1(\Omega) \geq \frac{(R_{m})^{n-2}}{(R_{M})^{n-1}} \frac{\left\lbrace 2+a-\sqrt{a^2+4\,a}\right\rbrace }{2\sqrt{a+1}},
\end{align*}
where $a:= \tan^{2} \alpha,$ $R_{m}:= \min  R(u),$ $R_{M}:= \max R(u)$ and $R(u):= \max \left\lbrace r : x \in \Omega, x = (r,u) \right\rbrace $, $(r,u)$ represents the polar coordinates on $\mathbb{R}^n.$ 
\end{thm}

In Theorem \ref{thm: star shaped domain}, we give a lower bound for the first non-zero steklov eigenvalue in a star-shaped domain $\Omega$ in $(\mathbb{S}^{n}, ds^{2}),$ following the idea of Kuttler and Sigillito. In particular, given a star-shaped domain $\Omega$ with respect to a point $p$ in $(\mathbb{S}^{n}, ds^{2})$, let $R_m$ be the minimum distance of the $\partial\Omega$ from $p$ and $B(R_m)$ be the ball of radius $R_m$ centered at $p$. We show that  $\mu_{1} (\Omega) \geq C \, \mu_{1}\left( B\left( R_{m}\right)\right)$ for some constant $C$ defined in Section 3.  \\

Let $(B_R, g)$ be a ball of radius $R$ in $\mathbb{R}^{n}$ with a rotationally invariant metric $g = dr^2 + \sigma^{2}(r) \, du^{2},$ where $du^{2}$ is the standard metric on $\mathbb{S}^{n-1}.$ 

Escobar \cite{E3} proved that the first non-constant eigenfunction for the Steklov problem on $B_R$ has the form
$\phi (r,u) = \psi(r) \, e(u),$
where $ e(u)$ and $\psi(r)$ satisfy 
\begin{align*}
&\Delta_{\mathbb{S}^{n-1}} e + (n-1) \, e =0, \\
&\frac{1}{\sigma^{n-1}(r)} \frac{d}{dr}\left(\sigma^{n-1}(r)\, \frac{d}{dr} \, \psi(r) \right) - \frac{(n-1)\, \psi(r) }{\sigma^{2}(r)} = 0 \quad \mbox{ in } (0,R), \\
&\psi'(R)= \mu_{1}(B_R) \, \psi(R), \quad \psi(0)=0.
\end{align*}
Using this result he proved that the first non-zero eigenvalue $\mu_{1}$ for two dimensional ball $B_R$ is $\mu_{1}= \frac{1}{\sigma(R)}$.

For $n$-dimensional ball $(B_{R}, g)$ with rotationally invariant metric $g=dr^2 + \sigma^{2}(r) \, du^{2}$, Montano \cite{M} proved the following.
\begin{itemize}
\item If $Ric(g) \geq 0$, then
$\mu_{1} \leq \left( \frac{R}{\sigma (R)}\right)^{n+1} \, \frac{1}{R}.$
\item If $Ric(g) \leq 0$, then
$\mu_{1} \geq \left( \frac{R}{\sigma (R)}\right)^{n+1} \, \frac{1}{R}.$
\end{itemize}

In Theorem \ref{thm: for bounded radial curvature}, we derive a two sided bound for $\mu_{1}$ on the ball with rotationally invariant metric and with bounded radial curvature.

\section{Preliminaries}

In this section we fix some notations which will be used in the proof of the main results.

\begin{defn}
Let $M$ be a Riemannian manifold and $p \in M$. The injectivity radius at $p$, denoted by $\operatorname{inj}(p)$, is defined as  $$\operatorname{inj}(p):= \sup_{r} \left\lbrace r>0 | \exp_{p} : B(0, r)\rightarrow  B(p, r) \mbox{ is a diffeomorphism}\right\rbrace. $$  
\end{defn}



Let $\Omega$ be a star-shaped bounded domain in $(\mathbb{S}^{n}, ds^{2})$ with respect to a point $p \in \Omega$ and $\Omega \subset B(p, \operatorname{inj}(p))$, a geodesic ball of radius $\operatorname{inj}(p)$. Let $\partial\Omega$ be the smooth boundary of $\Omega$ with unit outward normal $\nu.$ Since $\Omega$ is star-shaped with respect to the point $p$ and with smooth boundary, for every point $q \in \partial\Omega,$ there exists unique unit vector $u \in T_{p}\mathbb{S}^{n}$ and $R_{u} >0$ such that $q = \exp_{p} (R_{u} \, u)$. In notation we have the following.
\begin{align*}
\partial\Omega &= \left\lbrace \exp_{p} (R_{u} \, u): u \in T_{p}\mathbb{S}^{n}, \Vert u \Vert =1 \right\rbrace , \\
\Omega &= \left\lbrace \exp_{p} (t \, u): u \in T_{p}\mathbb{S}^{n}, \Vert u \Vert =1, 0 \leq t < R_{u} \right\rbrace. 
\end{align*}
In geodesic normal coordinates, $\Omega$ and $\partial\Omega$ can be written as
\begin{align*}
\partial\Omega &= \left\lbrace ( R_{u}, u)  : u \in T_{p} \mathbb{S}^{n}, \Vert u \Vert=1\right\rbrace,\\
\Omega\backslash \left\lbrace p\right\rbrace &= \left\lbrace (r, u)  : u \in T_{p} \mathbb{S}^{n}, \Vert u \Vert=1,  0 < r < R_{u}\right\rbrace.
\end{align*}
Let $R_{m}:=\mbox{ min } {R_{u}}, \, R_{M}:=\mbox{ max } {R_{u}}$.

For a Riemannian manifold $(M,g)$ with constant sectional curvature $k$, we know that
\begin{align*}
g =
\begin{cases}
dr^{2} + r^{2} \, g_{\mathbb{S}^{n-1}} & \text{ for } k=0, \\
dr^{2} + \frac{1}{k} \sin ^{2}(\sqrt{k}\, r) \, g_{\mathbb{S}^{n-1}} & \text{ for } k>0, \\
dr^{2} + \frac{1}{|k|}  \sinh ^{2}(\sqrt{|k|}\, r) \, g_{\mathbb{S}^{n-1}} & \text{ for } k<0,
\end{cases}
\end{align*}
with respect to the geodesic normal coordinates about a point, where $g_{\mathbb{S}^{n-1}}$ is the standard metric on the $(n-1)$-dimensional unit sphere and $ r \leq \frac{\pi}{\sqrt{k}}$ for $k >0.$ 
Then for any smooth function $f$ defined on $\overline{\Omega}$, 
$$\|\nabla f\|^2 = \left(\frac{\partial f}{\partial r} \right)^2 + \frac{1}{\sin^{2}r }\|\overline{\nabla} f\|^2,$$
where $\overline{\nabla}f$ represents the component of ${\nabla}f$ tangential to $\mathbb{S}^{n-1}$, the $(n-1)$-dimensional unit sphere. \\

Let $\partial_r$ be the radial vector field starting at $p,$ the centre of $\Omega$ and $\nu$ be the unit outward normal to $\partial\Omega$. Since $\Omega$ is a star-shaped bounded domain, for any point $q \in \partial\Omega$,  $\cos (\theta(q)) = \langle \nu(q), \partial_r(q) \rangle > 0 $. Therefore $\theta(q) < \frac{\pi}{2}$ for all $q \in \partial\Omega$. By compactness of $\partial\Omega$ there exist $\alpha$ such that $0 \leq \theta(q) \leq \alpha  < \frac{\pi}{2}$ for all $q \in \partial\Omega$. \\

Let $(B_R,g)$, $R>0$ be a ball in $\mathbb{R}^n$ centered at a point $p$, with a rotationally invariant metric
\begin{align*}
g = dr^2 + \sigma^{2}(r) \, du^{2},
\end{align*}
where $du^{2}$ is the standard metric on $\mathbb{S}^{n-1}$ and $ \sigma : [0, R) \rightarrow [0, \infty)$ is a non-zero smooth function satisfying $\sigma(0)=0, \,\sigma'(0)=1$ and $\sigma(r)>0, r \in (0, R]$. 

Let $q \in (B_R,g)$ and $X_{q} \perp \nabla r(q)$ be a unit tangent vector. The sectional curvature $K (\nabla r, X)$ is called the radial curvature and we know that
\begin{align*}
K (\nabla r, X) = - \frac{\sigma''(r)}{\sigma (r)}.
\end{align*}
Since $K (\nabla r, X)$ depends only on $r,$ we will denote it by $K(r).$

\section{ Main Results}

We begin this section by stating the main results.

The following theorem gives a sharp lower bound for the first non-zero steklov eigenvalue for a star-shaped domain in $(\mathbb{S}^{n}, ds^{2})$.

\begin{thm} \label{thm: star shaped domain}
Let $\Omega \subset \mathbb{S}^{n}, \nu$ and $\alpha$ are as in Section 2. Let $a = \max \left\lbrace \frac{\|\overline{\nabla} R_{u}\|^2}{\sin^{2}(R_{u})}\right\rbrace $. Then the first eigenvalue $\mu_{1}(\Omega)$ of \eqref{Steklov problem} satisfies the following inequality
\begin{align} \label{main ineq: lower bound}
\mu_{1} (\Omega) &\geq \left(\frac{R_{m}}{R_{M}} \right) \left( \frac{(2+a) - \sqrt{a^{2}+4\,a}}{2}\right)  \frac{\sin^{n-1}\left(R_{m}\right)}{\sec (\alpha) \, \sin^{n-1}\left(R_{M}\right)} \, \mu_{1}\left( B\left( R_{m}\right)\right),
\end{align}
where $B\left( R_{m}\right)$ is the geodesic ball of radius $R_{m}$ centered at $p$. Furthermore, equality occurs if and only if $\Omega$ is a geodesic ball.
\end{thm}

Next theorem gives a sharp lower bound as well as a sharp upper bound for the first non-zero steklov eigenvalue on the ball with rotationally invariant metric and with bounded radial curvature.

\begin{thm} \label{thm: for bounded radial curvature}
 Let $(B_R,g)$,  $K(r)$ are as in Section 2 and $(B_R,\operatorname{can}_{k})$ be a ball in $\mathbb{R}^n$ of radius $R$, with canonical metric $\operatorname{can}_{k}=dr^2 + {\sin_{k}}^{2}(r) \, du^{2}$. Then the following holds.
\begin{enumerate}[\rm(i)]
\item
If $K(r) \geq k$. Then
 \begin{align} \label{main ineq: for bounded radial curvature}
\mu_1(B_R,\operatorname{can}_{k}) \leq \mu_1(B_R,g)\leq \left( \frac{\sin_{k}R}{\sigma(R)}\right)^{n+1} \, \mu_1(B_R,\operatorname{can}_{k}).
 \end{align}

\item If $K(r) \leq k$. Then
\begin{align*}
\left( \frac{\sin_{k}R}{\sigma(R)}\right)^{n+1} \mu_1(B_R,\operatorname{can}_{k})\leq \mu_{1}(B_R,g)\leq \mu_{1} (B_R,\operatorname{can}_{k}).
\end{align*}

\end{enumerate} 
 
 Furthermore, equality holds if and only if $(B_R,g)$ is isometric to the $(B_R,\operatorname{can}_{k}).$
\end{thm}  
 
The following lemma is crucial to prove the Theorem \ref{thm: star shaped domain}.

\begin{lemma} \label{lem:inequality for sin}
For $x \in [0, \frac{\pi}{2}]$,
\begin{enumerate}[\rm(i)]
\item  $\sin ax \geq a\,\sin x $ \mbox{ for } $0 \leq a \leq 1$,
\item  $\sin ax \leq a\,\sin x $ \mbox{ for } $a \geq 1$.
\end{enumerate}
\end{lemma} \vspace{0.1 cm}

\begin{proof}[Proof]~ Define $g(x) = \sin ax - a\,\sin x $. Observe that $g(0)=0$ for all $a \in \mathbb{R}$.
	
\begin{enumerate}[\rm(i)]
\item We show that $g$ is an increasing function on $[0, \frac{\pi}{2}]$, for $0 \leq a \leq 1$. As $0 \leq ax \leq x \leq \frac{\pi}{2}$ and $\cos ax \geq \cos x$, it follows that $g'(x) \geq 0$. \\

\item Let $ a \geq 1$. Then $ ax \geq x $ for $x \in [0, \frac{\pi}{2}]$. We divide the proof into two cases.
 
 \begin{itemize}
 \item For $ax \leq \frac{\pi}{2}$, proof follows by using the same argument we used for $0 \leq a \leq 1.$ 
 
 \item For $ax \geq \frac{\pi}{2}$, we will show that $a\, \sin x \geq 1 $. Note that $a\, \sin x$ is an increasing function for $0 \leq x \leq \frac{\pi}{2}$. So it is enough to prove that $a\, \sin \left(\frac{\pi}{2\,a} \right)  \geq 1 $. Since $\frac{1}{a} \leq 1$, we have $ \sin \left(\frac{\pi}{2\,a} \right) \geq \frac{1}{a} \sin \left( \frac{\pi}{2}\right)= \frac{1}{a}$. Therefore $\sin ax \leq 1 \leq a\, \sin \left(\frac{\pi}{2\,a} \right) \leq a\,\sin x $. 
\end{itemize}

\end{enumerate}
\end{proof}
\begin{rmk}
If $\sin ax = a\,\sin x $ for all $x \in [0, \frac{\pi}{2}]$, then it can be easily seen that either $a=0$ or $a=1$.
\end{rmk}


\begin{proof}[\textbf{Proof of Theorem \ref{thm: star shaped domain}}]~
By variational characterization of $\mu_{1}$, we have
\begin{align*}
\mu_{1}(\Omega)= \inf \left\lbrace \frac{\int_\Omega{\|\nabla u\|^2}\, dv}{\int_{\partial\Omega}{u^2}\, ds} : \int_{\partial\Omega}{u}\, ds =0, \, u(\neq 0) \in C^1(\overline{\Omega})\right\rbrace.
\end{align*}
For a continuously differential real valued function $f$ defined on $\overline{\Omega}$, we find a lower bound for $ \int_\Omega{\|\nabla f\|^2}\, dv $ and an upper bound for $  \int_{\partial\Omega}{f^2}\, ds $ to find a lower bound for  $ \frac{\int_\Omega{\|\nabla f\|^2}\, dv}{\int_{\partial\Omega}{f^2}\, ds} $.

\break  

\begin{itemize}
\item
\textbf{A Lower bound for $\int_\Omega{\|\nabla f\|^2}\, dv$.}
\end{itemize}
Let $f$ be a continuously differential real valued function defined on $\overline{\Omega}$. Then for $q \in \Omega$, 
$$\|\nabla f\|^2 = \left(\frac{\partial f}{\partial r} \right)^2 + \frac{1}{\sin^{2}r }\|\overline{\nabla} f\|^2.$$
Therefore
\begin{align*}
\int_\Omega{\|\nabla f\|^2}\, dv = \int_{U_{p}\Omega} \int_{0}^{R_{u}} \left[ \left(\frac{\partial f}{\partial r} \right)^2 + \frac{1}{\sin^{2}r }\|\overline{\nabla} f\|^2\right] \sin^{n-1}r \, dr\, du.
\end{align*}
Let $r = \frac{\rho \, R_{u}}{R_{m}}$. Then the integral is written as
\begin{align*}
\int_\Omega{\|\nabla f\|^2}\, dv &= \int_{U_{p}\Omega} \int_{0}^{R_{m}} \left[\left\lbrace \left(\frac{R_{m}}{R_{u}} \right)^{2} +\|\overline{\nabla} R_{u}\|^2\left(\frac{\rho}{R_{u}\, \sin \left(\frac{\rho \, R_{u}}{R_{m}} \right) } \right)^2 \right\rbrace  \left(\frac{\partial f}{\partial \rho} \right)^2 + \frac{1}{\sin^{2}\left(\frac{\rho \, R_{u}}{R_{m}} \right) }\|\overline{\nabla} f\|^2 \right. \\
& \qquad \left. - \frac{2\, \rho}{R_{u}\, \sin^{2}\left(\frac{\rho \, R_{u}}{R_{m}} \right)} \frac{\partial f}{\partial \rho}\, \langle\overline{\nabla} f, \overline{\nabla} R_{u}\rangle\right] \sin^{n-1}\left(\frac{\rho \, R_{u}}{R_{m}} \right) \, \left( \frac{R_{u}}{R_{m}}\right) \, d\rho\, du \\
& =\int_{U_{p}\Omega} \int_{0}^{R_{m}} \left[\left\lbrace \left(\frac{R_{m}}{R_{u}} \right) + \frac{\|\overline{\nabla} R_{u}\|^2}{R_{u}\,R_{m} }\left(\frac{\rho}{\sin \left(\frac{\rho \, R_{u}}{R_{m}} \right) } \right)^2 \right\rbrace  \left(\frac{\partial f}{\partial \rho} \right)^2 + \frac{R_{u}}{R_{m}\,\sin^{2}\left(\frac{\rho \, R_{u}}{R_{m}} \right) }\|\overline{\nabla} f\|^2 \right. \\
& \qquad \left. - \frac{2\, \rho}{R_{m}\, \sin^{2}\left(\frac{\rho \, R_{u}}{R_{m}} \right)} \frac{\partial f}{\partial \rho}\, \langle\overline{\nabla} f, \overline{\nabla} R_{u}\rangle\right] \sin^{n-1}\left(\frac{\rho \, R_{u}}{R_{m}} \right) \, d\rho\, du. 
\end{align*}
Next we estimate $\langle\overline{\nabla} f, \overline{\nabla} R_{u}\rangle$. For any function $\beta^2$ on $\overline{\Omega}$, Cauchy-Schwartz inequality gives
\begin{align*}
- \frac{2\, \rho}{R_{m}\, \sin^{2}\left(\frac{\rho \, R_{u}}{R_{m}} \right)} \left( \frac{\partial f}{\partial \rho}\right) \, \langle\overline{\nabla} f, \overline{\nabla} R_{u}\rangle & \geq - \frac{1}{\beta^2}\frac{\|\overline{\nabla} R_{u}\|^2}{R_{u}\,R_{m} }\left(\frac{\rho}{\sin \left(\frac{\rho \, R_{u}}{R_{m}} \right) } \right)^2 \left(\frac{\partial f}{\partial \rho} \right)^2 -\frac{\beta^2 \, R_{u}}{R_{m}\,\sin^{2}\left(\frac{\rho \, R_{u}}{R_{m}} \right) }\|\overline{\nabla} f\|^2.
\end{align*}
Thus
\begin{align} \nonumber
\int_\Omega{\|\nabla f\|^2}\, dv & \geq \int_{U_{p}\Omega} \int_{0}^{R_{m}} \left[\left\lbrace \left(\frac{R_{m}}{R_{u}} \right) - \left(\frac{1}{\beta^2} - 1\right)  \frac{\|\overline{\nabla} R_{u}\|^2}{R\,R_{m} }\left(\frac{\rho}{\sin \left(\frac{\rho \, R_{u}}{R_{m}} \right) } \right)^2 \right\rbrace  \left(\frac{\partial f}{\partial \rho} \right)^2 \right. \\\label{grad f expression}
& \qquad \left. + \frac{R_{u}\, \left( 1- \beta^2 \right) }{R_{m}\,\sin^{2}\left(\frac{\rho \, R_{u}}{R_{m}} \right) }\|\overline{\nabla} f\|^2  \right] \sin^{n-1}\left(\frac{\rho \, R_{u}}{R_{m}} \right) \, d\rho\, du. 
\end{align}
Note that $0 \leq \frac{\rho}{R_{m}} \leq 1 \leq \frac{R_{u}}{R_{m}}$ and $0 \leq \rho,R_{u} \leq \frac{\pi}{2}$. Hence it follow from Lemma \ref{lem:inequality for sin} that
\begin{align*}
\frac{\rho}{R_{m}} \sin(R_{u})  \leq \sin\left( \frac{\rho\, R_{u}}{R_{m}}\right) \leq \frac{R_{u}}{R_{m}} \sin(\rho).
\end{align*} 
This gives 
\begin{align} \label{sin inequality1}
\frac{\sin(\rho)}{\sin\left( \frac{\rho\, R_{u}}{R_{m}}\right)} \geq \frac{R_{m}}{R_{u}} \qquad \mbox{ and } \qquad - \left(\frac{\rho}{\sin \left(\frac{\rho \, R_{u}}{R_{m}} \right) } \right)^2 \geq - \left( \frac{R_{m}}{\sin(R_{u})}\right)^2, \qquad \mbox{  for } \rho >0.
\end{align}
Since $0 \leq \rho \leq \frac{\rho\, R_{u}}{R_{m}} \leq R_{u} \leq \frac{\pi}{2}$, we have
 \begin{align} \label{sin inequality2}
0 \leq \sin^{n-1}(\rho) \leq \sin^{n-1}(\frac{\rho\, R_{u}}{R_{m}})\leq 1.
\end{align}
We assume $\beta^{2}<1$ and by substituting above inequalities in \eqref{grad f expression}, we get
\begin{align} \nonumber
\int_\Omega{\|\nabla f\|^2}\, dv & \geq \int_{U_{p}\Omega} \int_{0}^{R_{m}} \left[\left\lbrace \left(\frac{R_{m}}{R_{u}} \right) - \left(\frac{1}{\beta^2} - 1 \right)  \frac{\|\overline{\nabla} R_{u}\|^2}{R_{u}\,R_{m} }\left(\frac{R_{m}}{\sin(R_{u})} \right)^2 \right\rbrace  \left(\frac{\partial f}{\partial \rho} \right)^2 \right. \\\nonumber
& \qquad \left. + \frac{R_{u}\, \left( 1- \beta^2 \right) }{R_{m}} \left(\frac{R_{m}}{R_{u}\,\sin(\rho)} \right)^2  \|\overline{\nabla} f\|^2  \right] \sin^{n-1}\left(\rho\right) \, d\rho\, du \\\nonumber
& \geq \int_{U_{p}\Omega} \int_{0}^{R_{m}} \left(\frac{R_{m}}{R_{u}} \right) \left[\left\lbrace 1 - \left( \frac{1}{\beta^2} - 1\right)  \frac{\|\overline{\nabla} R_{u}\|^2}{\sin^{2}(R_{u}) } \, \right\rbrace  \left(\frac{\partial f}{\partial \rho} \right)^2 \right. \\\nonumber
& \qquad \left. + \frac{ \left( 1- \beta^2 \right) }{\sin^{2}(\rho)} \|\overline{\nabla} f\|^2  \right] \sin^{n-1}\left(\rho\right) \, d\rho\, du. 
\end{align}
Recall that $a = \max \left\lbrace \frac{\|\overline{\nabla} R_{u}\|^2}{\sin^{2}(R_{u})}\right\rbrace $. Hence
\begin{align*}
\int_\Omega{\|\nabla f\|^2}\, dv & \geq \left(\frac{R_{m}}{R_{M}} \right)\int_{U_{p}\Omega} \int_{0}^{R_{m}} \left[\left\lbrace 1 - \left(\frac{1}{\beta^2} - 1\right)  a \, \right\rbrace  \left(\frac{\partial f}{\partial \rho} \right)^2 \right. \\
& \qquad \left. + \frac{ \left( 1- \beta^2 \right) }{\sin^{2}(\rho)} \|\overline{\nabla} f\|^2  \right] \sin^{n-1}\left(\rho\right) \, d\rho\, du. 
\end{align*}
By solving the equation $ 1 - \left(\frac{1}{\beta^2} - 1\right) a =  1- \beta^2 $ for $\beta^2$ we see that
$$ \beta^{2} = \frac{-a +\sqrt{a^{2}+4\,a}}{2}.$$
Therefore $$1 - \left(\frac{1}{\beta^2} - 1\right) a =1- \beta^{2} = \frac{(2+a) - \sqrt{a^{2}+4\,a}}{2} > 0. $$
From this it follows that
\begin{align} \nonumber
\int_\Omega{\|\nabla f\|^2}\, dv & \geq \left(\frac{R_{m}}{R_{M}} \right) \left( \frac{(2+a) - \sqrt{a^{2}+4\,a}}{2}\right)  \int_{U_{p}\Omega} \int_{0}^{R_{m}} \left[\left(\frac{\partial f}{\partial \rho} \right)^2 \right. \\ \nonumber
& \qquad \left. + \frac{1}{\sin^{2}(\rho)} \|\overline{\nabla} f\|^2  \right] \sin^{n-1}\left(\rho\right) \, d\rho\, du\\ \label{grad f inequality}
&= \left(\frac{R_{m}}{R_{M}} \right) \left( \frac{(2+a) - \sqrt{a^{2}+4\,a}}{2}\right)  \int_{{B(R_{m})}}  {\|\nabla f\|^2}\, dv.
\end{align} \vspace{0.1 cm}  \\
\begin{itemize}
\item
\textbf{An upper bound for $\int_{\partial\Omega}{f^2}\, ds$.}
\end{itemize}
\vspace{0.1 cm}   
Recall that the Riemannian volume measure on $\partial\Omega$, denoted $ds$, is given by $ds = \sec (\theta)\, \sin^{n-1}\left(R_{u}\right)\, du$ (see \cite{TF}), where $\theta$ and $R_{u}$ are defined in Section 2. Then
\begin{align*}
\int_{\partial\Omega}{f^2}\, ds = \int_{U_{p}\Omega} f^2 \, \sec (\theta)\, \sin^{n-1}\left(R_{u}\right)\, du. 
\end{align*}
By using the fact that $\sin^{n-1} R_{m} \leq \sin^{n-1} R_{u} \leq \sin^{n-1} R_{M} $ and substituting $r = \frac{\rho \, R_{u}}{R_{m}}$, this integral becomes
\begin{align} \nonumber
\int_{\partial\Omega}{f^2}\, ds &\leq \frac{\sec (\alpha) \, \sin^{n-1}\left(R_{M}\right)}{\sin^{n-1}\left(R_{m}\right)}\int_{U_{p}\Omega} f^2 \, \, \sin^{n-1}\left(R_{m}\right)\, du \\ 
\label{f^2 inequality}
& = \frac{\sec (\alpha) \, \sin^{n-1}\left(R_{M}\right)}{\sin^{n-1}\left(R_{m}\right)}\int_{S(R_{m})} f^2 \, ds.
\end{align}
If $f$ satisfies $\int_{S(R_{m})} f \, ds = 0$, by inequalities \eqref{grad f inequality} and \eqref{f^2 inequality} we have
\begin{align} \nonumber
\frac{\int_\Omega{\|\nabla f\|^2}\, dv}{\int_{\partial\Omega}{f^2}\, ds} &\geq \left(\frac{R_{m}}{R_{M}} \right) \left( \frac{(2+a) - \sqrt{a^{2}+4\,a}}{2}\right)  \frac{\sin^{n-1}\left(R_{m}\right)}{\sec (\alpha) \, \sin^{n-1}\left(R_{M}\right)}\frac{\int_{{B(R_{m})}}  {\|\nabla f\|^2}\, dv}{\int_{S(R_{m})} f^2 \, ds} \\ \label{inequality for vari char for f}
&\geq \left(\frac{R_{m}}{R_{M}} \right) \left( \frac{(2+a) - \sqrt{a^{2}+4\,a}}{2}\right)  \frac{\sin^{n-1}\left(R_{m}\right)}{\sec (\alpha) \, \sin^{n-1}\left(R_{M}\right)} \mu_{1}\left( B\left( R_{m}\right)\right).
\end{align}
Next we find a test function for the Steklov problem on $B(R_{m})$ in terms of a first eigenfunction for the Steklov problem on $\Omega$ and use the above inequality to get the desired bound.

Let $\psi_1$ be a first eigenfunction for the Steklov problem \eqref{Steklov problem} on $\Omega$. Define 
\begin{align}
\psi := \psi_1 - \overline{\psi_1} \qquad \mbox{ where } \qquad \overline{\psi_1}= \frac{\int_{S(R_{m})} \psi_1 \, ds}{Vol(S(R_{m}))}.
\end{align} 
Observe that $\int_{S(R_{m})} \psi \, ds = 0$. Therefore from \eqref{inequality for vari char for f}, we have
\begin{align} \label{inequality for vari char for eigen fun}
\frac{\int_\Omega{\|\nabla \psi\|^2}\, dv}{\int_{\partial\Omega}{\psi^2}\, ds}
&\geq \left(\frac{R_{m}}{R_{M}} \right) \left( \frac{(2+a) - \sqrt{a^{2}+4\,a}}{2}\right)  \frac{\sin^{n-1}\left(R_{m}\right)}{\sec (\alpha) \, \sin^{n-1}\left(R_{M}\right)} \mu_{1}\left( B\left( R_{m}\right)\right).
\end{align}
Next by using the fact that $\int_{\partial\Omega} \psi_1 \, ds = 0$ we obtain
\begin{align} \nonumber
\frac{\int_\Omega{\|\nabla \psi\|^2}\, dv}{\int_{\partial\Omega}{\psi^2}\, ds} &= \frac{\int_\Omega{\|\nabla \psi_1\|^2}\, dv}{\int_{\partial\Omega}\left( \psi_{1}^{2} + \overline{\psi_{1}}^{2}\right) \, ds} \\\nonumber
&\leq \frac{\int_\Omega{\|\nabla \psi_1\|^2}\, dv}{\int_{\partial\Omega} {\psi_{1}^{2}}  \, ds} \\ \label{var char inequality for eigenval.}
& = \mu_{1} (\Omega).
\end{align}
Combining the inequality \eqref{inequality for vari char for eigen fun} and \eqref{var char inequality for eigenval.} we obtain
\begin{align*} 
\mu_{1} (\Omega) &\geq \left(\frac{R_{m}}{R_{M}} \right) \left( \frac{(2+a) - \sqrt{a^{2}+4\,a}}{2}\right)  \frac{\sin^{n-1}\left(R_{m}\right)}{\sec (\alpha) \, \sin^{n-1}\left(R_{M}\right)} \, \mu_{1}\left( B\left( R_{m}\right)\right).
\end{align*}
If $\Omega$ is a geodesic ball of radius $R$, then $R_{m}= R_{M}=R $ and $a = \alpha = 0$, hence equality holds in \eqref{main ineq: lower bound}. 
Next assume that equality hold in \eqref{main ineq: lower bound}, then equality holds in \eqref{sin inequality1} and \eqref{sin inequality2}. This implies $R_u{}= R_{m}.$ Hence $\Omega$ is a geodesic ball.
This proves the theorem.
\end{proof} 

The following lemma is used to prove the Theorem \ref{thm: for bounded radial curvature}.
 
Let
\begin{align*}
\sin_{k}r =
\begin{cases}
r & \text{ for } k=0, \\
\frac{\sin{\sqrt{k}}r}{\sqrt{k}} & \text{ for } k>0, \\
\frac{\sinh{\sqrt{-k}}r}{\sqrt{-k}} & \text{ for } k<0.
\end{cases}
\end{align*}

\begin{lemma} \label{lem: lemma for bounded radial curvature}
Let $(B_R,g)$, $R>0$ be a ball in $\mathbb{R}^n$ with a rotationally invariant metric $g = dr^2 + \sigma^{2}(r) \, du^{2}$ and $K(r)$ be the radial sectional curvature such that $K(r) \geq k$ for some $k \in \mathbb{R}$. We assume $ R \leq \frac{\pi}{\sqrt{k}}$ when $k >0.$ Then the following holds.

\begin{enumerate}[\rm(i)]
\item $\sigma(r)\leq \sin_{k} (r)$ for $r\in [0,R]$,
\item  $\displaystyle \lim_{r\rightarrow 0}\frac{\sigma(r)}{\sin_{k} (r)}=1$ and
\item $\frac{\sigma(r)}{\sin_{k} (r)}$ is a decreasing function.
\end{enumerate}
\end{lemma} \vspace{0.1 cm}

\begin{proof}[Proof]~ \begin{enumerate}[\rm(i)]
\item Let $k \in \mathbb{R}$ and observe that $\sin_{k}r$ is the solution of
\begin{align*}
y''(r) + k \, y(r) = 0, \, y(0)=0, \, y'(0)=1, \, r\in [0,R].
\end{align*}
Since $K(r) \geq k$ and 
\begin{align*}
\sigma ''(r) + K(r) \, \sigma(r)= 0, \, \sigma(0)=0, \, \sigma'(0)=1,\, r\in [0,R],
\end{align*}
it follows from Sturm Comparison Theorem \cite{DC} that  $\sigma(r)\leq \sin_{k} (r)$ for $r\in [0,R]$.  \\

\item Since $\sigma'(0)=1=\sin_{k}'(0)$, it follows from L'Hospital  Rule that
\begin{align*}
\lim_{r\rightarrow 0}\frac{\sigma(r)}{\sin_{k} (r)} = \lim_{r\rightarrow 0}\frac{\sigma'(r)}{\sin_{k}' (r)}=1.
\end{align*}  

\item For $k<0,$
$$
\left( \frac{\sigma(r)}{\sinh{\sqrt{-k}}r}\right)'= \frac{\sigma'(r) \, \sinh{\sqrt{-k}}r - {\sqrt{-k}} \, \sigma(r)\, \cosh{\sqrt{-k}}r}{\sinh^{2}{\sqrt{-k}}r}
$$
and
\begin{align*}
&\left( \sigma'(r) \, \sinh{\sqrt{-k}}r - {\sqrt{-k}} \, \sigma(r)\, \cosh{\sqrt{-k}}r\right)' \\ &= \sigma''(r) \, \sinh{\sqrt{-k}}r + {\sqrt{-k}} \, \sigma'(r) \, \cosh{\sqrt{-k}}r - {\sqrt{-k}} \, \sigma'(r)\, \cosh{\sqrt{-k}}r + k \, \sigma(r)\, \sinh{\sqrt{-k}}r \\
&= (\sigma''(r) + k\,\sigma(r))\, \sinh{\sqrt{-k}}r \\
&\leq 0.
\end{align*}
Therefore $\sigma'(r) \, \sinh{\sqrt{-k}}r - {\sqrt{-k}} \, \sigma(r)\, \cosh{\sqrt{-k}}r$ is a decreasing function of $r$ and 
$$
\sigma'(r) \, \sinh{\sqrt{-k}}r - {\sqrt{-k}} \, \sigma(r)\, \cosh{\sqrt{-k}}r \leq 0 \mbox{ for } r \geq 0.
$$
Hence $\left( \frac{\sigma(r)}{\sinh{\sqrt{-k}}r}\right)' \leq 0.$ This implies that $\frac{\sigma(r)}{\sinh{\sqrt{-k}}r}$ is a decreasing function. \\

A similar computation as above proves that $\frac{\sigma(r)}{\sin_{k} (r)}$ is a decreasing function for $k \geq 0$ as well.
\end{enumerate}
\end{proof}

\begin{proof}[\textbf{Proof of Theorem \ref{thm: for bounded radial curvature}}]~ First we prove the theorem for $K(r) \geq k$.\\
Let $f$ be a continuously differential real valued function defined on $\overline{B}_R$. Then
\begin{align*}
{\|\nabla f\|^2} =  \left(\frac{\partial f}{\partial r} \right)^2 + \frac{1}{\sigma^{2}(r) }\|\overline{\nabla} f\|^2,
\end{align*}
where $\overline{\nabla} f$ is the component of $\nabla f$ tangential to $\mathbb{S}^{n-1}.$
Therefore
\begin{align*}
 \int_{B_R}{\|\nabla f\|^2}\, dv_g = \int_{B_R}  \left[ \left(\frac{\partial f}{\partial r} \right)^2 + \frac{1}{\sigma^{2}(r) }\|\overline{\nabla} f\|^2\right] \sigma^{n-1}(r) \, dr\, du. 
 \end{align*}
By Lemma \ref{lem: lemma for bounded radial curvature}, $\frac{\sigma(R)}{\sin_{k} (R)} \leq \frac{\sigma(r)}{\sin_{k} (r)} \leq 1$ for $0 < r \leq R,$ we have 
 \begin{align*}
\int_{B_R}{\|\nabla f\|^2}\, dv_g &\geq \int_{B_R}  \left[ \left(\frac{\partial f}{\partial r} \right)^2 + \left( \frac{\sigma(r)}{\sin_{k}r}\right)^{2} \, \frac{1}{\sigma^{2}(r)}\|\overline{\nabla} f\|^2\right] \left( \frac{\sigma(R)}{\sin_{k}R}\right)^{n-1} {\sin_k}^{n-1}r \, dr\, du \\
&= \left( \frac{\sigma(R)}{\sin_{k}R}\right)^{n-1}  \int_{B_R}  \left[ \left(\frac{\partial f}{\partial r} \right)^2 +  \frac{1}{{\sin_{k}}^{2}r}\|\overline{\nabla} f\|^2\right]  {\sin_k}^{n-1}r \, dr\, du \\
&= \left( \frac{\sigma(R)}{\sin_{k}R}\right)^{n-1}  \int_{B_R}{\|\nabla f\|^2}\, dv_{\operatorname{can}_{k}}, 
\end{align*}
and
\begin{align*}
\int_{\partial B_R} f^2\, ds_g &= \int_{\partial B_R} f^2\, \sigma^{n-1}(R) du\\
&= \left( \frac{\sigma(R)}{\sin_{k}R}\right)^{n-1} \int_{\partial B_R} f^2\, {\sin_k}^{n-1}(R) du\\ 
&= \left( \frac{\sigma(R)}{\sin_{k}R}\right)^{n-1} \int_{\partial B_R} f^2\, ds_{\operatorname{can}_{k}}.
\end{align*}
Hence
\begin{align} \label{inequality1}
\frac{\int_{B_R}{\|\nabla f\|^2}\, dv_g}{\int_{\partial B_R} f^2\, ds_g} \geq \frac{\int_{B_R}{\|\nabla f\|^2}\, dv_{\operatorname{can}_{k}} }{\int_{\partial B_R} f^2\, ds_{\operatorname{can}_{k}}}.
\end{align}
 Now we obtain a test function for the Steklov problem on $(B_R,{\operatorname{can}_{k}})$ in terms of a first eigenfunction for the Steklov problem on $(B_R,g)$.
 
Let $\psi_1$ be the first eigenfunction for the Steklov problem \eqref{Steklov problem} on $(B_R,g)$. Then $\int_{\partial B_R} \psi_1 \, ds_g = 0$. Define 
\begin{align*}
\psi := \psi_1 - \overline{\psi_1} \qquad \mbox{ where } \qquad \overline{\psi_1}= \frac{\int_{\partial B_R} \psi_1 \, ds_{\operatorname{can}_{k}}}{\int_{\partial B_R}\, ds_{\operatorname{can}_{k}}}.
\end{align*} 
Observe that $\int_{\partial B_R} \psi \, ds_{\operatorname{can}_{k}} = 0$. From \eqref{variational characterization} and \eqref{inequality1} it follows that
\begin{align}\nonumber
\mu_{1} (B_R,{\operatorname{can}_{k}}) &\leq \frac{\int_{B_R}{\|\nabla \psi\|^2}\, dv_{\operatorname{can}_{k}} }{\int_{\partial B_R} \psi^2\, ds_{\operatorname{can}_{k}}} \\ \nonumber
&\leq \frac{\int_{B_R}{\|\nabla \psi\|^2}\, dv_g }{\int_{\partial B_R} \psi^2\, ds_g} \\ \nonumber
&\leq  \frac{\int_{B_R}{\|\nabla \psi_{1}\|^2}\, dv_g }{\int_{\partial B_R}  \psi_{1}^{2} \, ds_g} \\ \label{ineq: bound1}
&= \mu_{1} (B_R,g).
\end{align}

Next we prove that $\mu_{1}(B_R,g)\leq \left( \frac{\sin_{k}R}{\sigma(R)}\right)^{n+1} \, \mu_1(B_R,{\operatorname{can}_{k}}).$

Note that
\begin{align}\label{eq.integ}
\int_{B_R}{\|\nabla f\|^2}\, dv_g &= \int_{B_R}  \left[ \left(\frac{\partial f}{\partial r} \right)^2 + \frac{1}{\sigma^{2}(r) }\|\overline{\nabla} f\|^2\right] \sigma^{n-1}(r) \, dr\, du. 
\end{align}
Lemma \ref{lem: lemma for bounded radial curvature} gives $\frac{\sigma(R)}{\sin_{k} (R)} \leq \frac{\sigma(r)}{\sin_{k} (r)} \leq 1$ for $0 < r \leq R.$ By substituting this in \eqref{eq.integ}, we get 
\begin{align*}
\int_{B_R}{\|\nabla f\|^2}\, dv_g &\leq \int_{B_R}  \left[ \left(\frac{\partial f}{\partial r} \right)^2 + \left( \frac{\sin_{k}r}{\sigma(r)}\right)^{2} \, \frac{1}{{\sin_{k}}^{2}(r)}\|\overline{\nabla} f\|^2\right] {\sin_k}^{n-1}r \, dr\, du \\
&\leq \int_{B_R}  \left[ \left(\frac{\partial f}{\partial r} \right)^2 + \left( \frac{\sin_{k}R}{\sigma(R)}\right)^{2} \, \frac{1}{{\sin_{k}}^{2}(r)}\|\overline{\nabla} f\|^2\right] {\sin_k}^{n-1}r \, dr\, du \\
&=\left( \frac{\sin_{k}R}{\sigma(R)}\right)^{2} \, \int_{B_R}  \left[\left( \frac{\sigma(R)}{\sin_{k}R}\right)^{2} \, \left(\frac{\partial f}{\partial r} \right)^2 +  \frac{1}{{\sin_{k}}^{2}(r)}\|\overline{\nabla} f\|^2\right] {\sin_k}^{n-1}r \, dr\, du. 
\end{align*}
Since $\frac{\sigma(R)}{\sin_{k} (R)} \leq 1,$ it follows that
\begin{align*}
\int_{B_R}{\|\nabla f\|^2}\, dv_g &\leq \left( \frac{\sin_{k}R}{\sigma(R)}\right)^{2} \, \int_{B_R}  \left[ \left(\frac{\partial f}{\partial r} \right)^2 +  \frac{1}{{\sin_{k}}^{2}(r)}\|\overline{\nabla} f\|^2\right] {\sin_k}^{n-1}r \, dr\, du \\
&= \left( \frac{\sin_{k}R}{\sigma(R)}\right)^{2} \, \int_{B_R}{\|\nabla f\|^2}\, dv_{\operatorname{can}_{k}}.
\end{align*}
Hence
\begin{align*}
\frac{\int_{B_R}{\|\nabla f\|^2}\, dv_g}{\int_{\partial B_R} f^2\, ds_g} \leq \left( \frac{\sin_{k}R}{\sigma(R)}\right)^{n+1} \,  \frac{\int_{B_R}{\|\nabla f\|^2}\, dv_{\operatorname{can}_{k}} }{\int_{\partial B_R} f^2\, ds_{\operatorname{can}_{k}}}.
\end{align*} 
By similar computation as in the previous case, we get
\begin{align} \label{ineq: bound2}
\mu_{1}(B_R,g)\leq \left( \frac{\sin_{k}R}{\sigma(R)}\right)^{n+1} \, \mu_1(B_R,{\operatorname{can}_{k}}).
\end{align} 
By combining \eqref{ineq: bound1} and \eqref{ineq: bound2}, we have
\begin{align*}
\mu_{1} (B_R,{\operatorname{can}_{k}}) \leq \mu_{1}(B_R,g)\leq \left( \frac{\sin_{k}R}{\sigma(R)}\right)^{n+1} \mu_1(B_R,{\operatorname{can}_{k}}).
\end{align*}
If $K_X \equiv k$, then by the uniqueness theorem of ordinary differential equation, $\sigma(r)= \sin_{k}r$ for $r \in [0,R]$ and equality holds for all inequalities in the proof. Hence
\begin{align*}
\mu_{1} (B_R,{\operatorname{can}_{k}}) = \mu_{1}(B_R,g) = \left( \frac{\sin_{k}R}{\sigma(R)}\right)^{n+1} \mu_1(B_R,{\operatorname{can}_{k}}).
\end{align*}
Next we assume that
\begin{align*}
\mu_{1} (B_R,{\operatorname{can}_{k}}) = \mu_{1}(B_R,g) = \left( \frac{\sin_{k}R}{\sigma(R)}\right)^{n+1} \mu_1(B_R,{\operatorname{can}_{k}}).
\end{align*}
Then  $\frac{\sigma(R)}{\sin_{k} (R)} = \frac{\sigma(r)}{\sin_{k} (r)} = \lim_{r\rightarrow 0} \frac{\sigma(r)}{\sin_{k} (r)} =1$ and $\overline{\psi_1} = 0.$ Therefore $\sigma(r)= \sin_{k}r$ for $r \in [0,R].$ Hence the theorem follows. \\

Proof for $K(r) \leq k$ follows from similar computation as the proof for the case $K(r) \geq k$.

\end{proof}

\section*{acknowledgment} 

I would like to thank Prof. G. Santhanam for the discussions and valuable comments that were helpful in carrying out this work.

\end{document}